\documentclass[11pt]{epiarticle}
\usepackage{epimath}

\setpapertype{A4}

\usepackage[english]{babel}

\title{Ramanujan-Fourier series of certain arithmetic functions of two variables}
\titlemark{Ramanujan-Fourier series of certain arithmetic functions of two variables}
\author{Noboru Ushiroya}
\authormark{N.\ Ushiroya}
\date{} 
\journal{}

\begin{document}

\maketitle




\begin{prelims}

\def\abstractname{Abstract}
\abstract{We study Ramanujan-Fourier series of certain arithmetic functions of two variables. We generalize Delange's theorem to the case of arithmetic functions of two variables and give sufficient conditions for pointwise convergence of Ramanujan-Fourier series of arithmetic functions of two variables. We also give several examples which are not obtained by trivial generalizations of results on Ramanujan-Fourier series of functions of one variable.}




\end{prelims}


\section{Introduction}

Let $c_q (n)$ denote the Ramanujan sums defined in \cite{Ramanujan} as
\begin{equation}
c_q(n)=\sum_{\substack{a=1 \\ (a,q)=1}}^q \exp (2 \pi i an/q) , \nonumber
\end{equation}
where  $q$ and $n$ are positive integers and $(a,q)$ is the greatest common divisor of $a$ and $q$. Ramanujan proved that $c_q(n)$ can be rewritten as 
\begin{equation}
\nonumber
c_q(n)=\sum_{d \mid (q,n)} \mu (q/d) d, 
\end{equation}
where $\mu$ is the M$\rm{\ddot{o}}$bius function. Hardy \cite{Hardy} proved that, for fixed $n$, $c_q(n)$ is a multiplicative function. In other words,
\begin{equation}
\nonumber
c_{q_1 q_2} (n)=c_{q_1} (n) c_{q_2} (n)
\end{equation}
holds for any $q_1,q_2 \in \mathbb{N}$ satisfying $(q_1,q_2)=1$.
 Let $f:\mathbb{N} \mapsto \mathbb{C} $ be an arithmetic function. Ramanujan \cite{Ramanujan} investigated its Ramanujan-Fourier series which is an infinite series of the form
\begin{equation}
f(n) \sim \sum_{q=1}^\infty a_q c_q (n),
\end{equation}
where $a_q$ are called the Ramanujan-Fourier coefficients of $f$, and he obtained the following results.
\begin{align}
& \frac{\sigma_s (n)}{n^s} =\zeta(s+1) \sum_{q=1}^\infty \frac{c_q (n)}{q^{s+1}} , \label{eq:sigma} \\
& \frac{\varphi (n)}{n} =\frac{1}{\zeta(2)} \sum_{q=1}^\infty \frac{\mu (q)}{\varphi_2 (q)} c_q (n) , \label{eq:varphi} \\
& \tau (n)=- \sum_{q=1}^\infty \frac{\log q}{q} c_q (n) , \label{eq:tau} \\
& r(n) =\pi \sum_{q=1}^\infty \frac{(-1)^{q-1}}{2q-1} c_{2q-1} (n) , \label{eq:r}
\end{align}
where $\sigma_s (n)=\sum_{d | n} d^s$ with $s>0$, \ $\zeta(s)$ is the Riemann zeta function, \ $\varphi(n)$ is Euler's totient function, \ $\varphi_s (n)= n^s \prod_{p | n} (1-1/p^s) $, \ $\tau (n)$ is the number of divisors of $n$ and $r(n)$ is the number of representations of $n$ as the sum of two squares.

Let $*$ denotes the Dirichlet convolution, that is, $(f*g)(n)=\sum_{d|n} f(d) g(n/d)$ for arithmetic functions $f,g$, and let $\omega(n)$ be the number of distinct prime divisors of $n$. 
We say that $f : \mathbb{N} \mapsto \mathbb{C}$ is a multiplicative function if $f$ satisfies $f(mn)=f(m ) f( n) $
for any $m, n \in \mathbb{N} $ satisfying $ (m,  n)=1.$

Delange \cite{Delange} proved the following theorem.
\begin{theorem}[\cite{Delange}]
\label{th:Delange}
Let $f(n)$ be an arithmetic function satisfying
\begin{equation}
\label{eq:Del-1}
\sum_{n=1}^\infty 2^{\omega(n)} \frac{|(f * \mu)(n)|}{n} < \infty.
\end{equation}
Then its Ramanujan-Fourier series is pointwise convergent and
\begin{equation}
\nonumber
f(n)=\sum_{q=1}^\infty a_q c_q (n)
\end{equation}
holds where 
\begin{equation}
\nonumber
a_q=\sum_{m=1}^\infty \frac{(f * \mu)(qm)}{qm} .  
\end{equation}
Moreover, if $f$ is a multiplicative function, then $a_q$ can be rewritten as
\begin{equation}
\label{eq:Del-2}
a_q=\prod_{p \in \mathcal{P}} \Bigl( \sum_{e=\nu_p(q)}^\infty \frac{(f * \mu)(p^e)}{p^e} \Bigr),
\end{equation}
where $\mathcal{P}$ is the set of prime numbers and $\nu_p(q)=\left\{ \begin{array}{lll} \alpha & \mathrm{if} & p^\alpha || n \\ 0 & \mathrm{if} & p \nmid n . \end{array}  \right.   $
\end{theorem}

Delange noted that, if $f$ is a multiplicative function, then the condition (\ref{eq:Del-1}) is equivalent to the condition: $\sum_{p \in \mathcal{P}} \sum_{e=1}^\infty |f(p^e)-f(p^{e-1})| / p^e < \infty $ since $(f* \mu) (p^e)=f(p^e)-f(p^{e-1})$ for $e \geqq 1$. Under this condition, we can directly calculate Ramanujan-Fourier coefficients $a_q$ for certain arithmetic functions by using (\ref{eq:Del-2}). For example, if we set $f(n)=\varphi (n)/n$, then we can easily calculate the right-hand side of (\ref{eq:Del-2}) and obtain $a_q=(\zeta (2))^{-1} \mu(q) / \varphi_2 (q)$ which coincide with the Ramanujan-Fourier coefficients of (\ref{eq:varphi}).

Many results concerning Ramanujan-Fourier series of arithmetic functions of one variable are obtained by many mathematicians hitherto, however, as for Ramanujan-Fourier series of arithmetic functions of two variables, to my knowledge, few results are known.
We would like to extend Delange's theorem to the case of functions of two variables and obtain several examples which are extensions of (\ref{eq:sigma}) $\sim$ (\ref{eq:r}).


\section{Some Results}

Let $f$, $g: \mathbb{N} \times \mathbb{N} \mapsto \mathbb{C}$ be arithmeric functions of two variables.
The Dirichlet convolution of $f$ and $g$ is defined as follows.
\begin{equation}
\nonumber
(f*g)(n_1,n_2)=\sum_{m_1 | n_1, \ m_2|n_2} f(m_1 , m_2) g( n_1/m_1,n_2/m_2 ). 
\end{equation}

We use the same notation $\mu$ for the function 
\begin{equation}
\nonumber
\mu(n_1, n_2)=\mu(n_1)  \mu(n_2),
\end{equation}
which is the inverse of the constant function $1$ under the Dirichlet convolution, that is, $(\mu * 1) (n_1, n_2)= \delta(n_1, n_2)  $ holds where $\delta (n_1, n_2)=1$ or $0 $ according to whether $n_1 =n_2=1$ or not.

We investigate Ramanujan-Fourier series of arithmetical functions of two variables along Delange's article (\cite{Delange}). We first establish the following theorem which is an extension of Theorem \ref{th:Delange} to the case of arithmetic functions of two variables.

\begin{theorem}
\label{th:th2-1}
Let $f(n_1,n_2)$ be an arithmetic function of two variables satisfying
\begin{equation}
\label{eq:th2-1-1}
\sum_{n_1,n_2=1}^\infty 2^{\omega(n_1)} 2^{\omega(n_2)} \frac{|(f*\mu) (n_1,n_2)|}{n_1 n_2} < \infty . 
\end{equation}
Then its Ramanujan-Fourier series is pointwise convergent and
\begin{equation}
\label{eq:th2-1-2}
f(n_1,n_2)=\sum_{q_1,q_2=1}^\infty a_{q_1,q_2} c_{q_1} (n_1)c_{q_2} (n_2) 
\end{equation}
holds where 
\begin{equation}
\label{eq:th2-1-a}
a_{q_1,q_2}  = \sum_{m_1,m_2=1}^\infty \frac{(f * \mu) (m_1 q_1, m_2 q_2) }{m_1 q_1 m_2 q_2}.  
\end{equation}
\end{theorem}

For the proof of the above theorem, we need the following lemma.

\begin{lemma}[\cite{Delange}]
\label{lem:Del}
For every positive integer $k$,     
\begin{equation}
\nonumber
\sum_{q \mid k} |c_q (n)| \leqq n 2^{\omega(k)}.    
\end{equation}
\end{lemma}

\begin{proof}[Proof of Theorem \ref{th:th2-1}.]

We proceed as in \cite{Delange}. We first note that (\ref{eq:th2-1-1}) implies the absolute convergence of the right-hand side of (\ref{eq:th2-1-a}). Next we show that the series $\sum_{q_1,q_2=1}^\infty a_{q_1,q_2} c_{q_1} (n_1)c_{q_2} (n_2) $ is absolutely convergent.
It is easy to see that
\begin{align*}
\sum_{q_1,q_2=1}^\infty |a_{q_1,q_2} c_{q_1} (n_1)c_{q_2} (n_2) |  & \leqq \sum_{q_1,q_2=1}^\infty \sum_{m_1,m_2=1}^\infty \Bigl| \frac{(f*\mu) (q_1 m_1, q_2 m_2)}{q_1 m_1 q_2 m_2}  c_{q_1} (n_1)c_{q_2} (n_2) \Bigr|  \\   
& =\sum_{k_1,k_2=1}^\infty W_{k_1,k_2} , 
\end{align*}
where 
\begin{align*}
W_{k_1,k_2} & = \sum_{\substack{m_1 q_1=k_1 \\ m_2 q_2=k_2 }} \Bigl| \frac{(f*\mu) (q_1 m_1, q_2 m_2)}{q_1 m_1 q_2 m_2} c_{q_1} (n_1)c_{q_2} (n_2) \Bigr|=  \Bigl|\frac{(f*\mu) (k_1, k_2)}{k_1 k_2} \Bigr| \sum_{\substack{q_1 \mid k_1 \\ q_2 \mid k_2}} |c_{q_1} (n_1) c_{q_2} (n_2) | . 
\end{align*}
By (\ref{eq:th2-1-1}) and Lemma \ref{lem:Del} we have
\begin{align*}
\sum_{k_1,k_2=1}^\infty W_{k_1,k_2} \leqq  \sum_{k_1,k_2=1}^\infty  \Bigl|\frac{(f*\mu) (k_1, k_2)}{k_1 k_2} \Bigr|  n_1 n_2 2^{\omega(k_1)} 2^{\omega(k_2)} \ll n_1 n_2 < \infty .
\end{align*}
Hence the series $\sum_{q_1,q_2=1}^\infty a_{q_1,q_2} c_{q_1} (n_1)c_{q_2} (n_2) $ is absolutely convergent. We have
\begin{equation}
\nonumber
\sum_{q_1,q_2=1}^\infty \sum_{m_1,m_2=1}^\infty \frac{(f*\mu) (q_1 m_1, q_2 m_2)}{q_1 m_1 q_2 m_2}  c_{q_1} (n_1)c_{q_2} (n_2) =\sum_{k_1,k_2=1}^\infty w_{k_1,k_2},
\end{equation}
where
\begin{equation}
\nonumber
w_{k_1,k_2} = \sum_{\substack{m_1 q_1=k_1 \\ m_2 q_2=k_2 }} \frac{(f*\mu) (q_1 m_1, q_2 m_2)}{q_1 m_1 q_2 m_2}  c_{q_1} (n_1)c_{q_2} (n_2)  = \frac{(f*\mu) (k_1, k_2)}{k_1 k_2} \sum_{\substack{q_1 \mid k_1 \\ q_2 \mid k_2}} c_{q_1} (n_1) c_{q_2} (n_2).
\end{equation}
Using the well known formula:
$ \sum_{q \mid k} c_{q} (n) = \varepsilon_{k} (n)$ 
where $\varepsilon_{k} (n)=\left\{ \begin{array}{lll} k & \mathrm{if} & k \mid n \\ 0 & \mathrm{if} & k \nmid n  \end{array}  \right. $
(\cite{Sivaramakrishnan}), we have
\begin{align*}
\sum_{k_1,k_2=1}^\infty w_{k_1,k_2} & =\sum_{k_1,k_2=1}^\infty  \frac{(f*\mu) (k_1, k_2)}{k_1 k_2} \varepsilon_{k_1} (n_1) \varepsilon_{k_2} (n_2) \\
& = \sum_{\substack{k_1 \mid n_1 \\ k_2 \mid n_2}} \frac{(f*\mu) (k_1, k_2)}{k_1 k_2} k_1 k_2=(f*\mu * 1) (n_1,n_2)=(f*\delta)(n_1,n_2)=f(n_1,n_2).
\end{align*}
Therefore (\ref{eq:th2-1-2}) holds.
This completes the proof of Theorem \ref{th:th2-1}.
\end{proof}

We say that $f : \mathbb{N} \times \mathbb{N} \mapsto \mathbb{C}$ is a multiplicative function of two variables if $f$ satisfies\begin{equation}
\nonumber
f(m_1 n_1, m_2  n_2)=f(m_1 , m_2 )  \ f( n_1,n_2) 
\end{equation}
for any $m_1 ,m_2,  n_1, n_2 \in \mathbb{N} $ satisfying $ (m_1 m_2, \  n_1 n_2)=1.$
It is well known that if $f$ and $g$ are multiplicative functions of two variables, then $f*g$ also becomes a multiplicative function of two variables. For an arithmetical function $f$ of two variable, the mean value $M(f)$ is defined by 
\begin{equation}
\nonumber
M(f)= \lim_{x,y \to \infty} \frac{1}{xy} \sum_{n_1 \leqq x, \ n_2 \leqq y} f(n_1,n_2) 
\end{equation}
if this limit exists.
Ushiroya \cite{Ushiroya1} proved the following theorem.
\begin{theorem}[\cite{Ushiroya1}]
\label{th:Ushi1}
Let $f$ be a multiplicative function of two variables satisfying
\begin{equation}
\nonumber
\sum_{p \in \mathcal{P} } \sum_{\substack{e_1 ,e_2  \geqq 0 \\ e_1+e_2 \geqq 1}} \frac{| (f * \mu) (p^{e_1},p^{e_2}) |}{ p^{e_1+e_2} }  < \infty .
\end{equation}
Then the mean value $M(f)$ exists and
\begin{equation}
\label{eq:th_u}
M(f)=\sum_{m_1,m_2=1}^\infty \frac{(f * \mu) (m_1 , m_2 ) }{m_1 m_2 }=\prod_{p \in \mathcal{P}}  \Bigl( \sum_{e_1,e_2=0}^\infty \frac{(f * \mu) (p^{e_1},p^{e_2}) }{ p^{e_1+e_2} }  \Bigr).
\end{equation}
\end{theorem}

We would like to investigate Ramanujan-Fourier series in the case when $f$ is a multiplicative function of two variables.
The following theorem is an extension of Theorem \ref{th:Delange} to the case of a multiplicative function of two variables.

\begin{theorem}
\label{th:th2-4}
Let $f$ be a multiplicative function of two variables satisfying
\begin{equation}
\label{eq:th2_4_ass1}
 \sum_{p \in \mathcal{P} } \sum_{\substack{e_1 ,e_2  \geqq 0 \\ e_1+e_2 \geqq 1}} \frac{|(f * \mu ) (p^{e_1}, p^{e_2})|}{p^{e_1+e_2}}<\infty.
\end{equation}
Then its Ramanujan-Fourier series is pointwise convergent and
\begin{equation}
\label{eq:th2-4-1}
f(n_1,n_2)=\sum_{q_1,q_2=1}^\infty a_{q_1,q_2} c_{q_1} (n_1)c_{q_2} (n_2) 
\end{equation}
holds where 
\begin{equation}
\label{eq:th2-4-a}
a_{q_1,q_2}  = \sum_{m_1,m_2=1}^\infty \frac{(f * \mu) (m_1 q_1, m_2 q_2) }{m_1 q_1 m_2 q_2} \\
\end{equation}
Moreover, if $M(f) \neq 0$ and $\{ q_1,q_2 \} >1$, where $\{q_1,q_2 \}$ denotes the least common multiple of $q_1$ and $q_2$
, then $a_{q_1,q_2} $ can be rewritten as follows.
\begin{align}
a_{q_1,q_2} & =\prod_{p \in \mathcal{P}} \Bigl( \sum_{e_1=\nu_p (q_1)} \sum_{e_2=\nu_p (q_2)} \frac{(f * \mu ) (p^{e_1}, p^{e_2})}{p^{e_1+e_2}} \Bigr) \nonumber \\
\label{eq:th2-4-2}
 & =M(f) \prod_{p \mid \{q_1,q_2 \}} \Bigl\{ \Bigl( \sum_{e_1=\nu_p (q_1)} \sum_{e_2=\nu_p (q_2)} \frac{(f * \mu ) (p^{e_1}, p^{e_2})}{p^{e_1+e_2}} \Bigr) / \Bigl( \sum_{e_1=0} \sum_{e_2=0} \frac{(f * \mu ) (p^{e_1}, p^{e_2})}{p^{e_1+e_2}} \Bigr) \Bigr\} ,
\end{align}
\end{theorem}

\begin{remark}
\label{rem:rem2-5}
By (\ref{eq:th_u}) and (\ref{eq:th2-4-a}), we have $a_{1,1}=M(f)$.
\end{remark}

\begin{proof}
First we note that, if $f(n_1,n_2)$ is a multiplicative function of two variables, then $f(n_1,n_2)=\prod_{p \in \mathcal{P}} f(p^{\nu_p (n_1)},p^{\nu_p (n_2)})$ holds.
Let $f$ satisfy (\ref{eq:th2_4_ass1}). Since $(n_1,n_2) \mapsto 2^{\omega(n_1)} 2^{\omega(n_2)} |(f*\mu) (n_1,n_2)| / n_1 n_2 $ is a multiplicative function of two variables, we have
\begin{align*}
& \sum_{n_1 \leqq x ,  \ n_2 \leqq y} 2^{\omega(n_1)} 2^{\omega(n_2)} \frac{|(f*\mu) (n_1,n_2)|}{n_1 n_2} \leqq \sum_{k,\ell \geqq 0} \prod_{p \in \mathcal{P}} 2^{\omega(p^{k})} 2^{\omega(p^\ell)} \frac{|(f*\mu) (p^k,p^\ell)|}{p^{k+\ell}} \\
& \ll \prod_{p \in \mathcal{P}} \sum_{k,\ell \geqq 0} \frac{|(f*\mu) (p^k,p^\ell)|}{p^{k+\ell}} = \prod_{p \in \mathcal{P}} \Bigl(1+\sum_{\substack{e_1 , \ e_2  \geqq 0 \\ e_1+e_2 \geqq 1}} \frac{|(f * \mu ) (p^{e_1}, p^{e_2})|}{p^{e_1+e_2}} \Bigr) \\
& \leqq \exp \Bigl( \sum_{p \in \mathcal{P} } \sum_{\substack{e_1 , \ e_2  \geqq 0 \\ e_1+e_2 \geqq 1}} \frac{|(f * \mu ) (p^{e_1}, p^{e_2})|}{p^{e_1+e_2}} \Bigr) < \infty
\end{align*}
where we have used the well-known inequality $1+x \leqq \exp (x)$ for $x\geqq 0$.
Since $f$ satisfies (\ref{eq:th2-1-1}), we see by Theorem \ref{th:th2-1} that the Ramanujan-Fourier series of $f$ is pointwise convergent and (\ref{eq:th2-4-1}) holds.

Next we prove (\ref{eq:th2-4-2}) under the condition $M(f) \neq 0$ and $\{ q_1,q_2 \} >1$.
For $i=1,2$, let $q_i=\prod_j p_j^{e_{ij}}$ $(e_{ij} \geqq 0)$ and $m_i=r_i \prod_j p_j^{d_{ij}} $ $(d_{ij} \geqq 0)$ be the prime factor decompositions of  $q_i$ and $m_i$ respectively where $p_j$'s are prime numbers satisfying  $ p_j \mid \{q_1, q_2 \}  $ and $r_i$'s are positive integers coprime to $ q_1 q_2$. 
Then we have
\begin{equation}
\nonumber
 a_{q_1,q_2}  = \sum_{m_1,m_2=1}^\infty \frac{(f * \mu) (m_1 q_1, m_2 q_2) }{m_1 q_1 m_2 q_2} =\sum_{\substack{d_{ij} \geqq 0 \\  r_i \geqq 1,  \ (r_i, \ q_1 q_2)=1}} \frac{(f* \mu) ( r_1 \prod_j p_j^{d_{1j}+e_{1j}},  \ \ r_2 \prod_j p_j^{d_{2j}+e_{2j}} ) }{r_1 r_2 \prod_j p_j^{d_{1j}+e_{1j}+d_{2j}+e_{2j}}} .
\end{equation}
Since $f*\mu$ is multiplicative, we have
\begin{align*}
a_{q_1,q_2}  & = \sum_{\substack{d_{ij} \geqq 0 \\  r_i \geqq 1,  \ (r_i, \ q_1 q_2)=1}} \frac{(f* \mu) (\prod_j p_j^{d_{1j}+e_{1j}}, \ \prod_j p^{d_{2j}+e_{2j}} ) }{\prod_j p_j^{d_{1j}+e_{1j}+d_{2j}+e_{2j}} } \times \frac{(f* \mu) ( r_1,  r_2) }{r_1 r_2} \\
 & = \Bigl( \sum_{d_{ij} \geqq 0 } \frac{(f* \mu) (\prod_j p_j^{d_{1j}+e_{1j}}, \ \prod_j p^{d_{2j}+e_{2j}} ) }{\prod_j p_j^{d_{1j}+e_{1j}+d_{2j}+e_{2j}} } \Bigr) \Bigl(\sum_{  r_i \geqq 1,  \ (r_i, \ q_1 q_2)=1} \frac{(f* \mu) ( r_1,  r_2) }{r_1 r_2} \Bigr) \\
 & =\prod_{p \mid \{q_1,q_2 \}} \Bigl( \sum_{e_1=\nu_p (q_1)} \sum_{e_2=\nu_p (q_2)} \frac{(f * \mu ) (p^{e_1}, p^{e_2})}{p^{e_1+e_2}} \Bigr)  \times \prod_{p \nmid q_1 q_2} \Bigl( \sum_{e_1=0} \sum_{e_2=0} \frac{(f * \mu ) (p^{e_1}, p^{e_2})}{p^{e_1+e_2}} \Bigr).
\end{align*}
Since the condition $M(f) \neq 0$ implies $\sum_{e_1, e_2=0} \frac{(f * \mu ) (p^{e_1}, p^{e_2})}{p^{e_1+e_2}} \neq 0$ for every $p \in \mathcal{P}$ by (\ref{eq:th_u}), we have
\begin{align*}
a_{q_1,q_2} & = \displaystyle{ \frac{\displaystyle{ \prod_{p \mid \{q_1,q_2 \}} \Bigl( \sum_{e_1=\nu_p (q_1)} \sum_{e_2=\nu_p (q_2)} \frac{(f * \mu ) (p^{e_1}, p^{e_2})}{p^{e_1+e_2}} \Bigr)  \prod_{p \in \mathcal{P}} \Bigl( \sum_{e_1=0} \sum_{e_2=0} \frac{(f * \mu ) (p^{e_1}, p^{e_2})}{p^{e_1+e_2}} \Bigr) } }{ \displaystyle{\prod_{p \mid \{q_1,q_2 \}} \Bigr( \sum_{e_1=0} \sum_{e_2=0} \frac{(f * \mu ) (p^{e_1}, p^{e_2})}{p^{e_1+e_2}} \Bigr)  } }  }  \\
 & =M(f) \prod_{p \mid \{q_1,q_2 \}} \Bigl\{ \Bigr( \sum_{e_1=\nu_p (q_1)} \sum_{e_2=\nu_p (q_2)} \frac{(f * \mu ) (p^{e_1}, p^{e_2})}{p^{e_1+e_2}} \Bigr) / \Bigl( \sum_{e_1=0} \sum_{e_2=0} \frac{(f * \mu ) (p^{e_1}, p^{e_2})}{p^{e_1+e_2}} \Bigr) \Bigr\}.
\end{align*}
This completes the proof of Theorem \ref{th:th2-4}.
\end{proof}


\section{Examples}
In this section, we give several examples. We begin with the following example.


\begin{example}
\label{ex:3-1}
Let $f(n_1,n_2)=\varphi (n_1 n_2 ) / n_1 n_2 $. Then its Ramanujan-Fourier series is pointwise convergent and
\begin{equation}
\nonumber
\frac{\varphi (n_1 n_2 )}{n_1 n_2 }=M(f) \sum_{q_1,q_2=1}^\infty \frac{\mu (q_1) \mu(q_2) }{\varphi( (q_1 , q_2 )) \tilde{\varphi} ( q_1 q_2 )}  c_{q_1} (n_1) c_{q_2} (n_2)
\end{equation}
holds where $M(f)=\prod_{p \in \mathcal{P} } (1-2 /p^2 +1 /p^3)$ and $\tilde{\varphi} (n) =\prod_{p \mid n} (p^2 +p-1).$  
\end{example}
\begin{proof}
First we show that $f$ satisfies (\ref{eq:th2_4_ass1}). If we set $g(n)=\varphi (n) /n=\prod_{p \mid n} (1-1/p)$, then $g(p^e)=1-1/p$ for $e \geqq 1$. From this we have for $e,e_1,e_2 \geqq 1$
\begin{eqnarray*}
& (f*\mu)(p^e,1) & =f(p^e,1)-f(p^{e-1},1)=g(p^e)-g(p^{e-1})=\left\{ \begin{array}{cll} -1/p & \mathrm{if} & e=1 \\ 0 & \mathrm{if} & e \geqq 2,  \end{array}  \right.  \\
& \ \ \ (f*\mu)(p^{e_1},p^{e_2}) & =f(p^{e_1},p^{e_2}) -f(p^{e_1-1},p^{e_2})-f(p^{e_1},p^{e_2-1})+f(p^{e_1-1},p^{e_2-1}) \\
&  & = g(p^{e_1+e_2})-2g(p^{e_1+e_2-1})+g(p^{e_1+e_2-2}) \\
&  & = \left\{ \begin{array}{cll} 1/p & \mathrm{if} & (e_1,e_2)=(1,1) \\ 0 & \mathrm{if} & e_1,e_2 \geqq1 \ \ \mathrm{and} \ \   e_1+e_2 \geqq 3. \end{array}  \right. 
\end{eqnarray*}
From the above and (\ref{eq:th_u}) we can calculate the mean value $M(f)$ as follows.
\begin{align*}
M(f)  & =\prod_{p \in \mathcal{P} } \Bigl(1+\sum_{e=1}^\infty \frac{(f*\mu)(p^e,1)+(f*\mu)(1,p^e)}{p^e} +\sum_{e_1,e_2=1}^\infty \frac{(f*\mu)(p^{e_1},p^{e_2})}{p^{e_1+e_2}}  \Bigr) \\
 & =\prod_{p \in \mathcal{P} } \Bigl(1-\frac{2}{p^2} +\frac{1}{p^3} \Bigr).
\end{align*}
Since $M(f) \neq 0$, we have by (\ref{eq:th2-4-2}) 
\begin{align*}
& a_{p,1}=M(f) \Bigl(-\frac{1}{p^2}+\frac{1}{p^3} \Bigr)/ \Bigl(1-\frac{2}{p^2}+\frac{1}{p^3} \Bigr)=- \frac{M(f)}{ p^2+p-1} , \\
& a_{p^k,1}=0 \ \ \mathrm{if} \ \ k \geqq 2, \\
& a_{p,p}=M(f) \Bigl(\frac{1}{p^3} \Bigr)/ \Bigl(1-\frac{2}{p^2}+\frac{1}{p^3} \Bigr) =\frac{M(f)}{(p-1) (p^2+p-1)}, \\
& a_{p^k,p^\ell}=0 \ \ \mathrm{if} \ \  k,\ell \geqq 1 \ \ \mathrm{and}  \ \ k+\ell \geqq 3 .
\end{align*}
Hence 
\begin{equation}
\label{eq:ex3-1-2}
a_{q_1,q_2}=M(f) \frac{\mu (q_1) \mu(q_2) }{\varphi( (q_1 , q_2 )) \tilde{\varphi} ( q_1 q_2 )} 
\end{equation}
holds if $ (q_1,q_2)=(p^k,p^\ell)$ where $k,\ell \geqq 0$. Since the function $(q_1,q_2) \mapsto a_{q_1,q_2} / M(f) $ is multiplicative, (\ref{eq:ex3-1-2}) holds for every $q_1,q_2 \in \mathbb{N}$. This completes the proof of Example \ref{ex:3-1}.
\end{proof} 

The formula (\ref{eq:varphi}) says that, if we set $f(n)=\varphi (n) / n $, then $f(n)=M(f) \sum_{q=1}^\infty (\mu (q) / \varphi_2 (q)) c_q (n)$ holds since $M(f)=6/\pi^2$ in this case. However, Example \ref{ex:3-1} shows that the function $f(n_1,n_2)=\varphi (n_1 n_2 ) / n_1 n_2 $ does not have a similar expression of Ramanujan-Fourier series to that of (\ref{eq:varphi}).
We want to find an arithmetic function $f$ which satisfies $f(n_1,n_2)=M(f) \sum_{q_1,q_2=1}^\infty \frac{\mu (q_1 q_2)}{\varphi_2 (q_1 q_2)}  c_{q_1} (n_1)c_{q_2} (n_2)$. In the following example, we give a function which satisfies the above relation. 

\begin{example}
\label{ex:3-2}
  Let $f(n_1,n_2)$ be the multiplicative function defined by     
\begin{equation}
\label{eq:ex3-2-1}
f(p^k,p^\ell)= \left\{ \begin{array}{lll} 1-p/(p^2+1) & \mathrm{if} & k \ \ \mathrm{or} \ \  \ell=0  \ \ \mathrm{and} \ \  k+\ell \geqq 1 \\ 1-2p/(p^2+1) & \mathrm{if} & k, \ell \geqq 1.  \end{array}  \right. 
\end{equation}
Then its Ramanujan-Fourier series is pointwise convergent and
\begin{equation}
\nonumber
f(n_1,n_2)=M(f) \sum_{q_1,q_2=1}^\infty \frac{\mu (q_1 q_2)}{\varphi_2 (q_1 q_2)} c_{q_1} (n_1)c_{q_2} (n_2)
\end{equation}
holds where
\begin{equation}
\nonumber
M(f)=\prod_{p \in \mathcal{P}} (1-\frac{2}{p^2+1}).
\end{equation}

\end{example}
\begin{proof}
By (\ref{eq:ex3-2-1}) it follows that
\begin{equation}
\label{eq:ex2-5-1}
(f * \mu)(p^k,p^\ell)=\left\{ \begin{array}{ll} -p/(p^2+1) & \mathrm{if}  \ \ (k,\ell)=(1,0) \ \ \mathrm{or} \ \  (0,1) \\ 0 & \mathrm{otherwise}.    \end{array}  \right. 
\end{equation}
From the above it is obvious that $f$ satisfies (\ref{eq:th2_4_ass1}). Hence its Ramanujan-Fourier series is pointwise convergent by Theorem \ref{th:th2-4}. Using (\ref{eq:th_u}) and (\ref{eq:ex2-5-1}) we have
\begin{equation}
\nonumber
M(f)=\prod_{p \in \mathcal{P}} \Bigl(1+\frac{(f*\mu)(p,1)+(f*\mu)(1,p)}{p} \Bigr) =\prod_{p \in \mathcal{P}} (1-\frac{2}{p^2+1}).
\end{equation}
As for the Ramanujan-Fourier coefficients of $f$, we have by (\ref{eq:th2-4-2}) and (\ref{eq:ex2-5-1})
\begin{align*}
& a_{p^k,1}  =a_{1,p^k}  =\left\{ \begin{array}{lll} M(f) (-\frac{1}{p^2+1})/( 1-\frac{2}{p^2+1})=  M(f) \frac{-1}{p^2-1}=M(f) \frac{\mu(p)}{\varphi_2 (p)} & \mathrm{if} & k=1 \\ 0=M(f) \frac{\mu(p^k)}{\varphi_2 (p^k)} & \mathrm{if} & k \geqq 2,  \end{array}  \right.  \\
& a_{p^k,p^\ell}  =0=M(f) \frac{\mu(p^k p^\ell)}{\varphi_2 (p^k p^\ell)}  \ \ \mathrm{if} \ \ k,\ell \geqq 1.
\end{align*}
Hence 
\begin{equation}
\label{eq:ex3-2-2}
a_{q_1,q_2}=M(f) \frac{\mu (q_1 q_2)}{\varphi_2 (q_1 q_2)}
\end{equation}
holds if $ (q_1,q_2)=(p^k,p^\ell)$ where $k,\ell \geqq 0$. Since the function $(q_1,q_2) \mapsto a_{q_1,q_2} / M(f) $ is multiplicative, (\ref{eq:ex3-2-2}) holds for every $q_1,q_2 \in \mathbb{N}$. This completes the proof of Example \ref{ex:3-2}.
\end{proof}

Next we consider the case in which $f(n_1,n_2)$ is defined by $f(n_1,n_2)=g((n_1,n_2))$
where $g:\mathbb{N} \mapsto \mathbb{C}$ is a multiplicative function.
Before describing Ramanujan-Fourier series of this case, we would like to cite the following theorem.
\begin{theorem}[\cite{Ushiroya1}]
Let $f(n_1,n_2)=g((n_1,n_2))$ where $g$ is a multiplicative function satisfying 
\begin{equation}
\label{eq:th3-1}
\sum_{p \in P} \sum_{e \geqq 1} \frac{|g(p^e)-g(p^{e-1})|}{p^{2e}} < \infty.
\end{equation}
Then the mean value $M(f)$ exists and
\begin{equation}
\label{eq:th3-3-M}
M(f)=\prod_{p \in \mathcal{P}} \Bigl(1+\sum_{e =1}^\infty \frac{g(p^e)-g(p^{e-1})}{p^{2e}}   \Bigr).
\end{equation}
\end{theorem}

\begin{remark}
\label{re:rem3-4}
We note that, if $f(n_1,n_2)=g((n_1,n_2))$, then it is easy to see that
\begin{equation}
\label{eq:rem3-4}
(f*\mu) (p^k,p^\ell)= \left\{ \begin{array}{ll} 
   g(p^k)-g(p^{k-1}) \ \ \ & \mathrm{if} \ \  k=\ell \geqq 1  \\ 0 &  \mathrm{otherwise}.   \end{array}  \right.
\end{equation}
Hence (\ref{eq:th3-1}) clearly implies (\ref{eq:th2_4_ass1}).
\end{remark}

Now we can show the following theorem concerning Ramanujan-Fourier series in the case $f(n_1,n_2)=g((n_1,n_2)).$ 
Let $k \vee \ell$ denote $\max(k,\ell).$

\begin{theorem}
\label{th:3-4}
Let $f(n_1,n_2)=g((n_1,n_2))$ where $g$ is a multiplicative function satisfying 
\begin{equation}
\label{eq:th3-2-1}
\sum_{p \in P} \sum_{e \geqq 1} \frac{|g(p^e)-g(p^{e-1})|}{p^{2e}} < \infty.
\end{equation}
Then its Ramanujan-Fourier series is pointwise convergent and
\begin{equation}
\nonumber
f(n_1,n_2)=\sum_{q_1,q_2=1}^\infty a_{q_1,q_2} c_{q_1} (n_1) c_{q_2} (n_2) 
\end{equation}
holds where $a_{1,1}=M(f)$ and 
\begin{equation}
\label{eq:th3-2-2}
a_{q_1, q_2}=M(f) \prod_{p \mid \{q_1, q_2 \}} \Bigl\{ \Bigl( \sum_{e=\nu_p (q_1) \vee \nu_p (q_2)}^\infty \frac{g(p^e)-g(p^{e-1})}{p^{2e}} \Bigr) / \Bigl(1+\sum_{e=1}^\infty \frac{g(p^e)-g(p^{e-1})}{p^{2e}} \Bigr) \Bigr\} 
\end{equation}
for $q_1,q_2$ such that $\{ q_1,q_2 \} >1.$
\end{theorem}

\begin{proof}
Remark \ref{re:rem3-4} says that (\ref{eq:th3-2-1}) implies (\ref{eq:th2_4_ass1}).
It is easy to see that (\ref{eq:th3-2-2}) holds by (\ref{eq:th2-4-2}) and (\ref{eq:rem3-4}). This completes the proof of Theorem \ref{th:3-4}.
\end{proof}

The following example is an extension of (\ref{eq:sigma}) to the case $f(n_1,n_2)=g((n_1,n_2)).$
\begin{example}
\label{ex:sigma_s}
Let $f(n_1,n_2)=\sigma_s ((n_1,n_2)) / (n_1,n_2)^s $ where $s>-1$. Then its Ramanujan-Fourier series is pointwise convergent and     
\begin{equation}
\nonumber
\frac{\sigma_s ((n_1,n_2))}{(n_1,n_2)^s}=\zeta(s+2) \sum_{q_1,q_2=1}^\infty \frac{c_{q_1} (n_1) c_{q_2} (n_2)}{\{q_1,q_2 \}^{s+2}} .
\end{equation}
\end{example}
\begin{proof}
If we set $g(n)=\sigma_s (n) / n^s$, then we see that
\begin{equation}
\nonumber
g(p^e)= \frac{1}{p^{es}} \frac{1-p^{(e+1)s}}{1-p^s}= \frac{1}{1-p^s} \Bigl(\frac{1}{p^{es}}-p^s \Bigr) \ \ \mathrm{for} \ \  e \geqq 0
\end{equation}
and  
\begin{equation}
\nonumber
g(p^e)-g(p^{e-1})= \frac{1}{1-p^s} \Bigl(\frac{1}{p^{es}}-\frac{1}{p^{(e-1)s}} \Bigr)=\frac{1}{1-p^s} \frac{1-p^s}{p^{es}} =\frac{1}{p^{es}}  \ \ \mathrm{ for } \ \  e \geqq 1.
\end{equation}
Therefore (\ref{eq:th3-1}) holds and we have by (\ref{eq:th3-3-M})
\begin{equation}
\nonumber
M(f)=\prod_{p \in \mathcal{P}} \Bigl(1+\sum_{e =1}^\infty \frac{1}{p^{(s+2)e}}   \Bigr)=\zeta(s+2).
\end{equation}
For $k$ and $\ell$ satisfying $k+\ell \geqq 1$ we also have by (\ref{eq:th3-2-2})
\begin{align*}
a_{p^k,p^\ell} & =M(f) \frac{\sum_{e=k \vee \ell}^\infty \frac{g(p^e)-g(p^{e-1})}{p^{2e}} }{1+\sum_{e=1}^\infty \frac{g(p^e)-g(p^{e-1})}{p^{2e}}}   = M(f) \frac{\sum_{e=k \vee \ell}^\infty \frac{1}{p^{(2+s)e}} }{1+\sum_{e=1}^\infty \frac{1}{p^{(2+s)e}}} \\
 & = M(f) \frac{ \frac{1}{p^{(2+s)(k \vee \ell)}} \frac{1}{1-\frac{1}{p^{2+s}}}}{\frac{1}{1-\frac{1}{p^{2+s}}}} = M(f) \frac{1}{p^{(2+s)(k \vee \ell)}}  = M(f) \frac{1}{\{ p^k,p^\ell \}^{2+s}} .
\end{align*}
Therefore $a_{q_1,q_2}=M(f) / \{ p^k,p^\ell \}^{2+s}  $ holds if $(q_1,q_2)=(p^k,p^\ell)$.
Since the function $(q_1,q_2) \mapsto a_{q_1,q_2}/M(f)$ is multiplicative, we obtain the desired result.
\end{proof}

If we set $s=0$ in Example \ref{ex:sigma_s}, then we obtain the following example.

\begin{example}
Let $f(n_1,n_2)=\tau ((n_1,n_2)) $ where $\tau (n)$ is the number of divisors of $n$. Then its Ramanujan-Fourier series is pointwise convergent and   
\begin{equation}
\nonumber
\tau ((n_1,n_2))= \zeta(2) \sum_{q_1,q_2=1}^\infty \frac{c_{q_1} (n_1) c_{q_2} (n_2)}{\{q_1,q_2 \}^2}. 
\end{equation}
\end{example}

The following example is an extension of (\ref{eq:varphi}) to the case $f(n_1,n_2)=g((n_1,n_2)).$
\begin{example}
Let $f(n_1,n_2)=\varphi ((n_1,n_2)) / (n_1,n_2) . $  Then its Ramanujan-Fourier series is pointwise convergent and  
\begin{equation}
\nonumber
\frac{\varphi ((n_1,n_2))}{(n_1,n_2)}=\frac{1}{\zeta(3)} \sum_{q_1,q_2=1}^\infty \frac{\mu( \{q_1,q_2 \})}{\varphi_3 (\{q_1,q_2 \})} c_{q_1} (n_1) c_{q_2} (n_2). 
\end{equation}
\end{example}
\begin{proof}
If we set $g(n)=\varphi (n) / n$, then $g(p^e)=1-1/p$ for $e \geqq 1$.
Hence (\ref{eq:th3-1}) holds and we have by (\ref{eq:th3-3-M})
\begin{equation}
\nonumber
M(f)=\prod_{p \in \mathcal{P}} \Bigl(1+\sum_{e =1}^\infty \frac{g(p^e)-g(p^{e-1})}{p^{2e}}   \Bigr)=\prod_{p \in \mathcal{P}} \Bigl(1+\frac{(1-1/p)-1}{p^{2}}   \Bigr)=\prod_{p \in \mathcal{P}} \Bigl(1-\frac{1}{p^3} \Bigr)=\frac{1}{\zeta(3)}.
\end{equation}
We also have  by (\ref{eq:th3-2-2})
\begin{equation}
\nonumber
a_{p^k,p^\ell} =M(f) \frac{\sum_{e=k \vee \ell}^\infty \frac{g(p^e)-g(p^{e-1})}{p^{2e}} }{1+\sum_{e=1}^\infty \frac{g(p^e)-g(p^{e-1})}{p^{2e}}}   = 0 \ \ \mathrm{if} \ \ k \vee \ell \geqq 2
\end{equation}
and
\begin{equation}
\nonumber
a_{p,1}=a_{1,p}=a_{p,p}= M(f) \frac{-1/p^3}{1-1/p^3}    = M(f) \frac{-1}{p^3(1-1/p^3)}.
\end{equation}
Therefore $a_{q_1,q_2}=M(f) \mu( \{q_1,q_2 \})/ \varphi_3 (\{q_1,q_2 \}) $ holds if $(q_1,q_2)=(p^k,p^\ell)$.
Since the function $(q_1,q_2) \mapsto a_{q_1,q_2}/M(f)$ is multiplicative, we obtain the desired result.

\end{proof}

The proof of the following example is similar to that of the previous example.
\begin{example}
\label{ex:varphi_s}
Let $f(n_1,n_2)=\varphi_s ((n_1,n_2)) / (n_1,n_2)^s  $ where $s>-1.$  Then its Ramanujan-Fourier series is pointwise convergent and  
\begin{equation}
\nonumber
\frac{\varphi_s ((n_1,n_2))}{(n_1,n_2)^s}=\frac{1}{\zeta(s+2)} \sum_{q_1,q_2=1}^\infty \frac{\mu( \{q_1,q_2 \})}{\varphi_{s+2} (\{q_1,q_2 \})} c_{q_1} (n_1) c_{q_2} (n_2). 
\end{equation}
\end{example}

If we set $s=0$ in Example \ref{ex:varphi_s}, then we obtain the following example.

\begin{example}
Let $f(n_1,n_2)=\delta ((n_1,n_2))=\begin{cases}
   1 \ \ \ if  \ \ (n_1,n_2)=1  \\ 0 \ \ \ if  \ \ (n_1,n_2)>1.  \end{cases}  $ Then its Ramanujan-Fourier series is pointwise convergent and
\begin{equation}
\nonumber
\delta((n_1,n_2))= \frac{1}{\zeta(2)} \sum_{q_1,q_2=1}^\infty \frac{\mu (\{q_1,q_2 \} )}{\varphi_2 (\{q_1,q_2 \} )} c_{q_1} (n_1) c_{q_2} (n_2) .
\end{equation}
\end{example}

The following example is an extension of (\ref{eq:r}) to the case $f(n_1,n_2)=g((n_1,n_2)).$
\begin{example}
\label{ex:ex3-10}
Let $f(n_1,n_2)=\frac{1}{4} r((n_1,n_2))$ where $r(n)=\# \{(A,B) \in \mathbb{Z} \times \mathbb{Z} ; A^2+B^2=n \} . $  
Then its Ramanujan-Fourier series is pointwise convergent and  
\begin{equation}
\nonumber
\frac{1}{4} r((n_1,n_2))=M(f) \sum_{q_1,q_2=1}^\infty \frac{\chi (\{ q_1,q_2 \})}{\{ q_1,q_2 \}^2} c_{q_1} (n_1) c_{q_2} (n_2) , \end{equation}
where $\chi (n)= \left\{ \begin{array}{cl} 
   0 & \mathrm{if} \ \ n \mathrm{ \ \ is \ \ even}  \\ (-1)^{\frac{n-1}{2}} &  \mathrm{if}  \ \ n \mathrm{ \ \ is \ \  odd}  \end{array}  \right.$ and \ $M(f)=\displaystyle{\prod_{p >2, \ p \in \mathcal{P}} \frac{1}{1-\chi (p) / p^2}} $.
\end{example}
\begin{proof}
Let $g(n)=\frac{1}{4} r(n)$ and $e \in \mathbb{N}$. Then $g(2^e)=1$ and
\begin{equation}
\nonumber
g(p^e)=
\left\{ \begin{array}{cll} e+1 & \rm{if} & p \equiv 1 \pmod{4} \\ 0 & \rm{if} & p \equiv 3 \pmod{4}  \ \ \mathrm{and} \ \ e \mathrm{ \ \ is \ \ odd} \\ 1 & \rm{if} & p \equiv 3 \pmod{4}  \ \ \mathrm{and} \ \ e \mathrm{ \ \ is \ \ even}. \end{array}  \right. 
\end{equation}
From this we have $g(2^e)-g(2^{e-1})=0$ and 
\begin{equation}
\nonumber
g(p^e)-g(p^{e-1})=\left\{ \begin{array}{rll} 1 & \rm{if} & p \equiv 1 \pmod{4} \\ -1 & \rm{if} & p \equiv 3 \pmod{4}  \ \ \mathrm{and} \ \ e \mathrm{ \ \ is \ \ odd} \\ 1 & \rm{if} & p \equiv 3 \pmod{4}  \ \ \mathrm{and} \ \ e \mathrm{ \ \ is \ \ even}. \end{array}  \right. 
\end{equation}
Hence (\ref{eq:th3-1}) clearly holds. Furthermore, we have for $k \geqq 1$
\begin{equation}
\label{eq:ex3-10-1}
\sum_{e=k }^\infty \frac{g(p^e)-g(p^{e-1})}{p^{2e}}=
\left\{ \begin{array}{lll} 0 & \mathrm{if} & p=2 \\ \sum_{e=k }^\infty \frac{1}{p^{2e}}=\frac{1}{p^{2k}} \frac{1}{1-1/p^2} & \mathrm{if} & p \equiv 1 \pmod{4} \\ \sum_{e=k }^\infty \frac{(-1)^e}{p^{2e}}=-\frac{1}{p^{2k}} \frac{1}{1+1/p^2} & \mathrm{if} & p \equiv 3 \pmod{4}  \ \ \mathrm{and} \ \ k \mathrm{ \ \ is \ \ odd} \\ \sum_{e=k }^\infty \frac{(-1)^e}{p^{2e}}=\frac{1}{p^{2k}} \frac{1}{1+1/p^2} & \mathrm{if} & p \equiv 3 \pmod{4}  \ \ \mathrm{and} \ \ k \mathrm{ \ \ is \ \ even}. \end{array}  \right. 
\end{equation}
Using (\ref{eq:th3-3-M}) and(\ref{eq:ex3-10-1}) we can calculate the mean value $M(f)$ as follows.
\begin{align*}
M(f) & = \prod_{p \equiv 1 (\bmod 4)} \Bigl(1+\frac{1}{p^2(1-1/p^2)} \Bigr) \prod_{p \equiv 3 (\bmod{4})} \Bigl(1-\frac{1}{p^2(1+1/p^2)} \Bigr) \\
 & = \prod_{p \equiv 1 (\bmod 4)} \frac{1}{1-1/p^2} \prod_{p \equiv 3 (\bmod 4)} \frac{1}{1+1/p^2}= \prod_{p >2, \ p \in \mathcal{P}} \frac{1}{1-\chi (p)/p^2} .
\end{align*}
Let $k,\ell$ be non-negative integers satisfying $k+\ell \geqq 1$.
Using (\ref{eq:th3-2-2}) and (\ref{eq:ex3-10-1}) we obtain
\begin{align*}
& a_{2^k,2^\ell}=0,  \\
& a_{p^k,p^\ell}=M(f) \frac{\frac{1}{p^{2(k \vee \ell)}} \frac{1}{1-1/p^2} }{ 1+\frac{1}{p^2-1} } = M(f) \frac{1}{p^{2(k \vee \ell)}}  \ \quad \mathrm{if} \ \ p \equiv 1 \pmod{4}, \\
& a_{p^k,p^\ell}=M(f) \frac{\frac{(-1)^{k \vee \ell}}{p^{2(k \vee \ell)}} \frac{1}{1+1/p^2} }{ 1-\frac{1}{p^2+1} } =M(f) \frac{(-1)^{k \vee \ell}}{p^{2(k \vee \ell)}}   \quad  \mathrm{if} \ \ p \equiv 3 \pmod{4}.
\end{align*}
Therefore $a_{q_1,q_2}=M(f) \chi (\{ q_1,q_2 \}) / \{ q_1,q_2 \}^2 $ holds if $(q_1,q_2)=(p^k,p^\ell)$ where $p \in \mathcal{P}$.
Since the function $(q_1,q_2) \mapsto a_{q_1,q_2}/M(f)$ is multiplicative, we obtain the desired result.
This completes the proof of Example \ref{ex:ex3-10}.

\end{proof}


\bibliographystyle{amsalpha}

\authoraddresses{
Noboru Ushiroya\\
National Institute of Technology, Wakayama College, \\
77 Noshima Nada Gobo Wakayama, Japan \\
\email ushiroya@wakayama-nct.ac.jp

}

\end{document}